\documentclass[11pt,a4paper,reqno]{amsart}
\usepackage[english]{babel}
\usepackage[applemac]{inputenc}
\usepackage[T1]{fontenc}
\usepackage{palatino}
\usepackage{verbatim}
\usepackage{amsmath}
\usepackage{amssymb}
\usepackage{amsthm}
\usepackage{amsfonts}
\usepackage{graphicx}
\usepackage{mathtools}
\usepackage{enumitem}

\usepackage[colorlinks = true, citecolor = black]{hyperref}
\author{Tuomas Orponen}
\title{Every locally finite Borel measure on $\R$ has conformal dimension zero}
\address{University of Helsinki, Department of Mathematics and Statistics}
\thanks{T.O. is supported by the Academy of Finland via the project Restricted families of projections and connections to Kakeya type problems, grant No. 274512.}
\email{tuomas.orponen@helsinki.fi}

\newcommand{\R}{\mathbb{R}}

\newcommand{\Z}{\mathbb{Z}}
\newcommand{\tn}{\mathbb{P}}

\newcommand{\calL}{\mathcal{L}}
\newcommand{\calD}{\mathcal{D}}
\newcommand{\calH}{\mathcal{H}}

\newcommand{\calG}{\mathcal{G}}

\newcommand{\spt}{\operatorname{spt}}
\newcommand{\Hd}{\dim_{\mathrm{H}}}

\newcommand{\E}{\mathbb{E}}

\newcommand{\card}{\operatorname{card}}

\newcommand{\Cd}{\dim_{\mathrm{C}}}

\numberwithin{equation}{section}

\theoremstyle{plain}
\newtheorem{thm}[equation]{Theorem}

\theoremstyle{definition}

\theoremstyle{remark}

\newcommand{\nref}[1]{(\hyperref[#1]{#1})}

\begin{document}

\begin{abstract} A result of P. Tukia from 1989 says that Lebesgue measure on $\R$ has conformal dimension zero: for every $\epsilon > 0$, there is a Borel set $G \subset \R$ of full Lebesgue measure, and a quasisymmetric homeomorphism $f \colon \R \to \R$ such that $\Hd f(G) < \epsilon$. In this short note, I show that the same is true for every locally finite Borel measure on $\R$. 
\end{abstract} 

\maketitle

\section{Introduction}

Let $(X,\mu,d)$ be a metric measure space. The \emph{conformal dimension} of $\mu$ is 
\begin{displaymath} \Cd \mu := \inf_{f} \Hd f_{\sharp}\mu, \end{displaymath}
where the $\inf$ is taken over all quasisymmetric homeomorphisms $f$ between $(\spt \mu,d)$ and any metric space $(Y,d')$. The notation $\Hd \nu$ stands for the (upper) Hausdorff dimension of $\nu$:
\begin{displaymath} \Hd \nu := \inf \{\Hd A : \nu(X \setminus A) = 0\}. \end{displaymath}
The concept of conformal dimension for measures is quite poorly understood: it currently not known, whether any measure has non-zero conformal dimension. This appears to be a hard problem for Lebesgue measure on $\R^{n}$, for $n \geq 2$. See the introduction of \cite{BO} for more information. 

For $n = 1$, the problem is not so hard. A construction of P. Tukia \cite{Tu} from 1989 shows that the Lebesgue measure on $[0,1]$, denoted by $\calL$, has conformal dimension zero (the restriction to $[0,1]$ is only for convenience, and proof works for Lebesgue measure on $\R$). The argument is the following. Pick $\epsilon > 0$ and find a doubling probability measure $\mu$ on $[0,1]$ with $\Hd \mu < \epsilon$ (such measures can be easily obtained via Riesz products, or the construction in \cite{GKS}). Define a quasisymmetric homeomorphism $f \colon [0,1] \to [0,1]$ by $f(t) := \mu([0,t])$, and note that $f_{\sharp}\mu = \calL$. Now $f^{-1} \colon [0,1] \to [0,1]$ is a quasisymmetric homeomorphism with $f^{-1}_{\sharp}(\calL) = \mu$. The proof is complete. 

What if $\calL$ is replaced by an arbitrary measure $\nu$ on $[0,1]$ or $\R$? The reasoning above shows that any pair of doubling probability measures on $[0,1]$ can be mapped to each other with a quasisymmetric self-homeomorphism of $[0,1]$ (start by mapping both measures to $\calL$), so the conformal dimension of any doubling measure on $[0,1]$ is zero. Also, if $\nu$ is an arbitrary measure on $\R$ with $\Hd \nu < 1$, a deep result of L. Kovalev \cite{Ko} implies that $\Cd \nu = 0$. So, the remaining problem concerns non-doubling measures $\nu$ with $\Hd \nu = 1$. While writing \cite{BO}, we believed that this would be a non-trivial problem, so we mentioned it explicitly after \cite[Question 1]{BO}. In the end, the solution turns out to be quite easy, and the purpose of this note is to record the proof.

\begin{thm}\label{main} Let $\nu$ be a locally finite Borel measure on $\R$. Then $\Cd \nu = 0$. \end{thm}

\subsection{Acknowledgements} We found the argument while discussing with Fredrik Ekstr\"om, but he did not wish to be an author. If the reader has any positive thoughts about this note, then they should equally be attributed to Fredrik.

\section{Proof of the main theorem}  

Fix $\epsilon > 0$, and any interval $I_{k} = [k,k + 1)$, $k \in \Z$. The plan is to construct increasing quasisymmetric homeomorphisms $f_{k} \colon I_{k} \to I_{k}$ such that $\Hd f_{k\sharp}(\nu|_{I_{k}}) < \epsilon$. This will be done so that $f \colon \R \to \R$, defined by $f|_{I_{k}} = f_{k}$, is also a quasisymmetric homeomorphism.

The first step will be to define, inductively, a doubling probability measure $\mu_{k} = \mu_{k,\epsilon}$ on $I_{k}$. Then, $f_{k}$ is defined by
\begin{equation}\label{f} f_{k}(t) := k + \mu_{k}([k,t]), \qquad t \in I_{k}. \end{equation}
Finally, it will be verified that $\calH^{\epsilon}(f_{k}(G_{k})) = 0$ for a Borel set $G_{k} \subset I_{k}$ with $\nu(I_{k} \setminus G_{k}) = 0$.

Let 
\begin{displaymath} \calD = \bigcup_{n \geq 0} \calD_{n} \end{displaymath}
be the collection of $4$-ary subintervals of $I_{k}$. Thus $\calD_{0} = \{I_{k}\}$, and $\calD_{n}$ is obtained from $\calD_{n - 1}$ by partitioning each interval in $\calD_{n - 1}$ into four half-open subintervals of length $4^{-n}$. Set $\mu_{k}(I_{k}) := 1$, and assume that $\mu_{k}(I)$ has already been defined for all intervals $I \in \calD_{n}$, for some $n \geq 0$. 

To continue, fix $I \in \calD_{n}$, and let $I^{1},\ldots,I^{4} \in \calD_{n + 1}$ be the children of $I$, from left to right. Let $\rho \in (0,\tfrac{1}{4})$ be a constant satisfying $4 \cdot \rho^{\epsilon/3} < 1$; note that $\rho$ is independent of $k$. Set 
\begin{equation}\label{form5} \mu_{k}(I^{1}) := \rho \cdot \mu_{k}(I) =: \mu_{k}(I^{4}). \end{equation}
If $\nu(I^{2}) \geq \nu(I^{3})$, set also $\mu_{k}(I^{2}) := \rho \cdot \mu_{k}(I)$; in the opposite case, set $\mu_{k}(I^{3}) := \rho \cdot \mu_{k}(I)$. The $\mu_{k}$-measure of the remaining interval, say $I^{j}$, needs to be $\mu_{k}(I^{j}) = [1 - 3\rho] \cdot \mu_{k}(I)$, so that
\begin{equation}\label{form1} \sum_{i = 1}^{4} \mu_{k}(I^{i}) = \mu_{k}(I). \end{equation}
holds. This completes the inductive definition of the set function $\mu_{k}$, and \eqref{form1} ensures that $\mu_{k}$ extends to a Borel probability measure on $I_{k}$. The point of the definition is that the intervals with most $\nu$ mass have the least $\mu_{k}$ mass.  

It is easy to verify that $\mu_{k}$ is a doubling measure, following the argument in \cite[Section 2.1]{GKS}. Here is the idea. Fix adjacent $4$-ary intervals $I^{1},I^{2} \in \calD_{n}$ of the same length; assume that $I^{1}$ is to the left from $I^{2}$. Let $J \in \calD_{k}$, $k < n$, be the smallest interval with $I^{1},I^{2} \subset J$. Then $I^{1} \subset J^{1}$ and $I^{2} \subset J^{2}$ for some adjacent, distinct $J^{1},J^{2} \in \calD_{k + 1}$ with $J^{1},J^{2}  \subset J$. After this stage, by adjacency, $I^{1}$ is always the right-most child of every descendant of $J^{1}$ between $I^{1}$ and $J^{1}$ (including $J^{1}$). The same holds for $I^{2},J^{2}$ with "right-most" replaced by "left-most". From the construction of $\mu_{k}$, it then follows that
\begin{equation}\label{doublingEstimate} \frac{\mu_{k}(I^{1})}{\mu_{k}(I^{2})} = \frac{\mu_{k}(J^{1})}{\mu_{k}(J^{2})} \in \left[\tfrac{\rho}{1 - 3\rho},\tfrac{1 - 3\rho}{\rho} \right]. \end{equation}
This proves that $\mu_{k}$ is doubling on $I_{k}$.  

Moreover, as long as each individual measure $\mu_{k}$, $k \in \Z$, is defined by the process above, respecting \eqref{form5} and with the constant $\rho$ fixed, then the sum $\mu = \sum \mu_{k}$ is a doubling measure on $\R$. To see this, consider two arbitrary adjacent $4$-ary intervals $I^{1},I^{2}$. If both are contained in a single $I_{k}$, then \eqref{doublingEstimate} holds for $\mu$. The situation $|I^{1}| = |I^{2}| \geq 1$ is also trivial, since each $\mu_{k}$ is a probability measure. So, the remaining case is, where $I^{1} \subset I_{k}$ and $I^{2} \subset I_{k + 1}$ for some $k \in \Z$. Then, with $4^{-n} := |I^{1}| = |I^{2}| \leq 1$, repeated application of condition \eqref{form5} (for both $\mu_{k}$ and $\mu_{k + 1}$) shows that 
\begin{displaymath} \mu(I^{1}) = \mu_{k}(I^{1}) = \rho^{n} = \mu_{k + 1}(I^{2}) = \mu(I^{2}). \end{displaymath}
This proves that $\mu$ is doubling on $\R$.

Fix $k \in \Z$. A Borel set $G_{k} \subset I_{k}$ is now constructed such that $\nu(I_{k} \setminus G_{k}) = 0$ and $\calH^{\epsilon}(f_{k}(G_{k})) = 0$. If $\nu(I_{k}) = 0$, take $G_{k} = I_{k}$. Otherwise, consider the probability space $(I_{k},\tn_{k})$, with $\tn_{k} = \nu/\nu(I_{k})$. Define the random variables $X^{n}_{k}$, $n \geq 1$, as follows. If $\omega \in I_{k}$, let $I_{n - 1} \in \calD_{n - 1}$ and $I_{n} \in \calD_{n}$ be the unique $4$-ary intervals containing $\omega$. Set
\begin{displaymath} X^{n}_{k}(\omega) := \begin{cases} +1, & \text{if } \mu(I_{n}) = \rho \cdot \mu(I_{n - 1}), \\ -1, & \text{if } \mu(I_{n}) = [1 - 3\rho] \cdot \mu(I_{n - 1}). \end{cases} \end{displaymath}
Write $S^{0}_{k} := 0$, and 
\begin{displaymath} S^{n}_{k} := \sum_{j = 1}^{n} X^{j}_{k}, \qquad n \geq 1. \end{displaymath}
Then $(S_{k}^{n})_{n \geq 0}$ is a sub-martingale. First of all,
\begin{equation}\label{form2} \E[S_{k}^{n} \mid S_{n - 1}] = \E[X^{n}_{k} \mid S_{n - 1}] + S_{n - 1}, \qquad n \geq 1. \end{equation} 
Second, note that $\sigma(S_{n - 1}) \subset \sigma(D_{n - 1})$, so 
\begin{equation}\label{form3} \E[X^{n}_{k} \mid S_{n - 1}] = \E[\E(X^{n}_{k} \mid \sigma(D_{n - 1})) \mid S_{n - 1}]. \end{equation}
So, to prove the sub-martingale property $\E[S_{k}^{n} \mid S_{n - 1}] \geq S_{n - 1}$, it remains by \eqref{form2}--\eqref{form3} to verify that $\E[X^{n}_{k} \mid \sigma(D_{n - 1})] \geq 0$. Fix $I \in D_{n - 1}$ with $\nu(I) > 0$, and let $I^{1},\ldots,I^{4} \in \calD_{n}$ be the children of $I$. Assume, for instance, that $\nu(I^{2}) \geq \nu(I^{3})$. Then,
\begin{displaymath} \E[X^{n}_{k} \mid I] = \frac{[\nu(I^{1}) + \nu(I^{2}) + \nu(I^{4})] - \nu(I^{3})}{\nu(I)} \geq \frac{\nu(I^{1}) + \nu(I^{4})}{\nu(I)} \geq 0. \end{displaymath}
If $\nu(I^{2}) < \nu(I^{3})$, then the same holds with the roles of $\nu(I^{2})$ and $\nu(I^{3})$ reversed. This proves that $(S_{k}^{n})_{n \geq 0}$ is a sub-martingale.  

The Azuma-Hoeffding inequality now says that
\begin{displaymath} \nu(\{S_{k}^{n} < -t\}) \leq e^{-t^{2}/(2n)}, \qquad t \geq 0. \end{displaymath}
In the current application, very crude estimates suffice: take $t = 2n^{3/4}$, and set $B_{k}^{n} := \{S_{k}^{n} < -2n^{3/4}\}$. Then $\nu(B^{n}) \leq \exp(-n^{1/2})$, so easily
\begin{equation}\label{B} \sum_{n = 0}^{\infty} \nu(B_{k}^{n}) < \infty. \end{equation}
Let $G_{k}^{n} := [0,1) \setminus B_{k}^{n}$, which is a union of certain $4$-ary intervals, say $\calG^{n}_{k} \subset \calD_{n}$. Fixing $I \in \calG^{n}_{k}$, $n \geq 1$, the variable $S_{k}^{n}(\omega)$ has constant value $S_{k}^{n}(I) \geq -n^{3/4}$ for $\omega \in I$. All the variables $X^{j}_{k}$ with $1 \leq j \leq n$ are also constant on $I$, and the constants are denoted by $X^{j}_{k}(I)$. Let $+(I) := \card \{1 \leq j \leq n : X^{j}_{k}(I) = +1\}$ and $-(I) := \card \{1 \leq j \leq n : X^{j}_{k}(I) = -1\}$, so that $+(I) - (-(I)) = S_{k}^{n}(I) \geq -2n^{3/4}$ and $+(I) + (-(I)) = n$. It follows that 
\begin{displaymath} +(I) \geq n/2 - n^{3/4}, \end{displaymath} 
and consequently
\begin{equation}\label{form4} \mu_{k}(I) = \rho^{+(I)} \cdot [1 - 3\rho]^{-(I)} \leq (\rho^{1/2 - n^{-1/4}})^{n} \leq \rho^{n/3} \end{equation}
for $n^{-1/4} < 1/6$. Since $|f_{k}(I)| = \mu_{k}(I)$, this implies that
\begin{equation}\label{form6} \calH^{\epsilon}_{\delta}(f_{k}(G_{k}^{n})) \leq \sum_{I \in \calG^{n}_{k}} \mu_{k}(I)^{\epsilon} \leq \card \calG^{n}_{k} \cdot \rho^{\epsilon n/3} \leq 4^{n} \cdot \rho^{\epsilon n/3} \end{equation}
assuming that $n$ is sufficiently large that \eqref{form4} holds, and $\rho^{n/3} < \delta$. Now, recall that $\rho^{\epsilon/3} < 1/4$, so the right hand side of \eqref{form6} tends to zero as $n \to \infty$. Consequently, setting
\begin{displaymath} G_{k} := \liminf_{n \to \infty} G^{n}_{k} = \bigcup_{m = 1}^{\infty} \bigcap_{n = m}^{\infty} G^{n}_{k}, \end{displaymath}
one finds that $\calH^{\epsilon}(f_{k}(G_{k})) = 0$, and \eqref{B} implies that $I_{k} \setminus G_{k} = \limsup_{n \to \infty} B^{n}_{k}$ satisfies $\nu(I_{k} \setminus G_{k}) = 0$. 

Since $\mu = \sum \mu_{k}$ is a doubling measure on $\R$, the map $f \colon \R \to \R$ defined by $f|_{I_{k}} := f_{k}$ is a quasisymmetric homeomorphism. Moreover, the Borel set 
\begin{displaymath} G := \bigcup_{k \in \Z} G_{k} \end{displaymath}
satisfies $\nu(\R \setminus G) = 0$ and $\calH^{\epsilon}(f(G)) = 0$. Hence $\Cd \nu < \epsilon$, and the proof of Theorem \ref{main} is complete.

\end{document}